\documentclass[12pt]{article}
\usepackage{amsfonts,amssymb,amsmath}

\newtheorem{theorem}{Theorem} 
\newtheorem{definition}{Definition}

\newtheorem{lemma}{Lemma}
 
\begin{document}
\title{Regularity of some class  of nonlinear transformations}

\maketitle

\noindent 
N. N. Ganikhodjaev \footnote{nasirgani@hotmail.com}\\
 Department of Mechanics and Mathematics,National University of Uzbekistan, Vuzgorodok, 700095, Tashkent, Uzbekistan and
 Centre for Computational and Theoretical Sciences,
Faculty of Science, International Islamic University Malaysia, 53100 Kuala Lumpur, Malaysia . 
\vskip 0.1cm
\noindent 
M. R. B. Wahiddin  \footnote{mridza@iiu.edu.my} \\
Centre for Computational and Theoretical Sciences,
Faculty of Science, International Islamic University Malaysia, 53100 Kuala Lumpur, Malaysia
\vskip 0.1cm
\noindent 
D.V.Zanin \footnote{zanin61@yandex.ru} \\
Department of Mechanics and Mathematics,National University of Uzbekistan, Vuzgorodok, 700095, Tashkent, Uzbekistan

\abstract{In this paper we consider quadratic stochastic operators designed on finite Abelian groups. It is proved that such operators have the property of regularity.}
\vskip 0.5mm
\noindent {\bf Mathematics Subject Classification }: 37A25 , 37N25 , 46T99 , 47H60.
\noindent {\bf Key words}:Quadratic stochastic operators; Finite Abelian group ; Regularity ; Ergodicity .

\section{Quadratic stochastic operators}
Let 
\begin{equation}
S^{n-1}=\{x=(x_1,\cdots,x_n)\in\mathbb{R}^n:\ \sum_{i=1}^n x_{i}=1\ ,\ x_i\geq0 \ \forall i=1,\cdots,n  \}
\end{equation}
be the $(n-1)$-dimensional simplex in $\mathbb{R}^n.$ The transformation
$V:S^{n-1}\rightarrow S^{n-1}$ is called a quadratic stochastic operator
(q.s.o.), if
\begin{equation}
(Vx)_k=\sum_{i,j=1}^np_{ij,k}x_ix_j
\end{equation}
where 
\begin{eqnarray} \nonumber
p_{ij,k} &\geq& 0, \\ \nonumber
p_{ij,k} &=& p_{ji,k}, \\
\sum_{k=1}^np_{ij,k}&=&1
\end{eqnarray}
for arbitrary $i,j,k \in\{1, \cdots ,n\}$. Such operators have applications in mathematical biology, namely theory
of heredity, where the coefficients $p_{ij,k}$ are interpreted as
coefficients of heredity [1-3].

Assume $\{V^kx:\ k=0,1,\cdots\}$ is the trajectory of the initial
point $x\in S^{n-1},$ where $V^{k+1}x=V(V^{k}x) $ for any 
$ k=0,1,\cdots $

\begin{definition}
A q.s.o. $V:S^{n-1}\rightarrow S^{n-1}$ is called
ergodic (respectively regular) if for any initial point $x\in S^{n-1}$
the limit 
$$\lim_{k\rightarrow\infty}\frac1k \sum_{i=0}^{k-1}V^ix$$
(respectively the limit $\lim\limits_{k\rightarrow\infty}V^{k}x$)
exists.
\end{definition}
Evidently, any regular q.s.o. $V$ has the ergodic property, but the
converse is not necessarily true.

To determine whether some q.s.o. is ergodic or regular is rather
complicated problem.

S. Ulam in [4] presupposed the assumption that any q.s.o. $V$ is
ergodic. Later M. Zakharevitch [5] showed, that this is an incorrect
hypothesis in general. More precisely, he proved that the q.s.o. $V,$
which is defined on the simplex
$$S^{2}=\{(x,y,z):\ x,y,z\geq0,\ x+y+z=1\}$$
by the formula

\[  V: \begin{cases}
\hat{x}=x^2+2xy;\\
\hat{y}=y^2+2yz;\\
\hat{z}=z^2+2xz.
\end{cases} \]
is not an ergodic q.s.o. and respectively is not regular.

Later in [6] necessary and sufficient conditions were established for
the ergodicity of the so-called Volterrian q.s.o.:
$$ V:   \begin{cases}
\hat{x}=x(1+ay-bz);\\
\hat{y}=y(1-ax+cz);\\
\hat{z}=z(1+bx-cy) 
\end{cases} $$
where $a,b,c\in[-1,1].$ For $a,b,c=1$ we have the initial example of
Zakharevitch.

\begin{theorem} Any q.s.o. in the above form is non-ergodic if and
only if the three parameters $a,b,c$ have the same sign.
\end{theorem}

Below we'll  construct one class of q.s.o. and prove that all such
q.s.o. are regular and respectively ergodic.

\section{The design of quadratic stochastic operators}

Let $G$ be a finite Abelian group and $S(G)$ be a set of all
probabilistic measures on $G.$ It is evident, that if $|G|=n,$  then
$S(G)$ coincides with $S^{n-1}  (1).$ 

Let further $H\subset G$ be a subgroup of $G$ and $\{g+H:\ g\in G\}$ be
the cosets of $H$ in $G.$ Assume $\mu\in S(G)$ is a fixed positive
measure, that is $\mu(g)>0$ for any $g\in G.$ Then we define the
coefficients $p_{fg,h},$ where $f,g,h\in G$ in the following way:
$$ p_{fg,h}=  \begin{cases}
\frac{\mu(g)}{\mu(f+g+H)},~ \mbox{if}~ \ h\in f+g+H;\\
0\,~ \mbox{otherwise}.
\end{cases} 
$$
It is easy to check that for arbitrary $f,g,h\in G$ the conditions $(3)$
are satisfied. It is also evident that if $H=\{e\},$ where $e$ is the neutral
element of group $G,$ then 
$$p_{fg,h}= \begin{cases}
1\ ~\mbox{if} \quad h=f+g;\\
0,~ \mbox{otherwise}
\end{cases} $$
and if $H=G,$ then 
$$p_{fg,h}=\mu(h) \quad \forall f,g\in G$$.
In the common case q.s.o. $V$ on $S(G)$ is defined as
$$(Vx)_{h}=\sum_{f,g\in G}p_{fg,h}x_fx_g$$
for all $h\in G,$ where $x=\{x_t,\ t\in G\}\in S(G)$ and $p_{fg,h}$ as
above.

If $H=\{e\}$ then the q.s.o. $V$ is defined as 
$$(Vx)_h=\sum_{f,g\in G,f+g=h}x_fx_g$$
and if $H=G,$ the q.s.o. is defined as
$$(Vx)_h=\mu(h)$$
for arbitrary $h\in G.$

Let us fix a positive measure $\mu\in S(G)$ and  subgroup $H$ of group $G.$ Assume $\mu_H$ is the factor-measure on factorgroup $G/H,$  that is 
$$\mu_{H}(g+H)=\sum_{h\in H}\mu(g+h)$$
for any $g\in G$ and $V_H$ is a q.s.o. on $S(G/H),$ which is defined by
measure $\mu_H.$ It is easy to show, that the trajectorial behaviour of
$V$ and $V_H$ are  similar,and  so it is enough to study the q.s.o.
generated by the  trivial subgroup.

Below we consider q.s.o. constructed by the trivial subgroup $H=\{e\}.$

\section{The main result}
Let $\nu\in S(G)$ be a Haar measure on $G $. As $G$ is a finite Abelian
group,the Haar measure $\nu $ on $G$ is a uniform distribution on $G$,that is 
it is a centre of the simplex  $S(G)$.
We prove the following 
\begin{theorem}
Almost all orbits tend to the center of the simplex.
\end{theorem}
{\it Proof:} The following lemma is a key ingredient:

\begin{lemma}
$||Vx||_{\infty}\leq ||x||_{\infty}.$
\end{lemma}
{\it Proof of lemma:} Let's define the following function
$$f(p)=f_n(p)=\max\limits_{\sum_0^nx_i=1,x_i\in[0,p]}\sum_0^n x_i^2 .$$
Notice  that there is no dependence on $n.$ Moreover, $f$  has an explicit form:
$$f(p)=kp^2+(1-kp)^2\ if\ p\in[\frac1{k+1},\frac1k].$$

We can prove this formula by induction .

The sum $\sum x_i^2$ cannot have a maximum in the interior of the domain (by the
Lagrange method). So, we can look for the maximum at the set
$x_n=p.$ It follows that 
$$f_n(p)=p^2+\max\limits_{x_i\in[0,p],\sum_0^{n-1}x_i=1-p}\sum_0^{n-1}
x_i^2$$
Let $x_i=(1-p)y_i.Then $ $\sum_0^{n-1}y_i=1$  and $y_i\in [0,\frac{p}{1-p}]$. So, the following recurrent formula is true:
$$f_n(p)=p^2+(1-p)^2f_{n-1}(\frac{p}{1-p}).$$
This allows an easy check of the expected formula.
Now it is easy to see that
$$f(p)-p=kp^2-p+(1-kp)^2=(1-kp)(1-kp-p)=(1-kp)(1-(k+1)p)\leq 0 ,$$
that is 
$$f(p)\leq p .$$

From Cauchy-Bunyakowski inequality:
$$(Vx)_i\leq \sqrt{\sum_{j\in G} x_j^2}\sqrt{\sum_{j\in G}
x_{i-j}^2}=\sum_{j\in G} x_j^2	$$

If $\max x_i\leq p$ then, $f(p)\leq p$ implies  $(Vx)_i\leq p$
or $\max (Vx)_i\leq p.$ So, the norm $||x||_{\infty}$ decreases on
the orbits. Lemma is thus  proven.

Actually, the set $\{x: ||x||_{\infty}\leq p\}$ is mapped to the set $\{x:
||x||_{\infty}\leq f(p)\}.$ We can easily check, that the iterations of the function
$f$ tend to its fixed points.
 
Now  consider a point $x.$ We consider the case, when a point in
orbit does not tend to a centre. It's $\omega-$limit set lies at a set
$\{x: g(x)=\frac1k\}.$ This $\omega-$limit set is finite. So, it is
nothing else, then a periodic orbit.

We consider the periodic orbits on the set $||x||_{\infty}=\frac1k.$

Since $||Vx||_{\infty}=\frac1k,$  there exists $k$ coordinates of $x$
equal to $\frac1k.$ All other coordinates are zeros. This is also true
for $Vx.$ So, there exists  sets $A=\{i:\ x_i=\frac1k\}$ and $B=\{j:\
(Vx)_j=\frac1k\},$ such that $|A|=|B|=k$ 

and 

\begin{enumerate}
\item $\forall i\in B$ $|(i-A)\bigcap A|=k $ $\rightarrow$ $i-A=A ,$
\item $\forall i\notin B$ $|(i-A)\bigcap A|=0 $ $\rightarrow$
$(i-A)\bigcap A=\emptyset.$
\end{enumerate}

From the second item: 

$$\forall i\notin B\ \forall x,y\in A:\ i-x\neq y ;$$
$$\forall i\notin B\ \forall j\in A+A ,\ i\neq j ;$$
$$A+A\subset B .$$

From the first item 
$$\forall i\in B\ \exists j,k \in A: \ i=j+k ;$$
$$B\subset A+A .$$

So, there exists  a set $A,$ such that $|A|=|A+A|=k,$ 
$$\forall i\in A+A:\ i-A=A .$$

{\bf Claim:} This is equivalent to the following property:
$$\forall i,j,k\in A : \quad i+j-k \in A .$$

{\bf Proof:} $\rightarrow$ $\forall x,y\in A$ $x+y-A=A.$ $\forall x,y,z\in
A\ x+y-z\in A.$

$\leftarrow$ $\forall x,y\in A$ $x+y-A\subset A.$ $\forall z\in A+A$ $z-A\subset A.$
Since $|z-A|=|A|,$ $z-A=A.$

$\forall a,i,j\in A$ $a+i-j\in A.$ $\forall i,j\in A$ $A+i-j\subset A.$
$$\forall i,j\in A\ A+i\subset A+j.$$
So, $\forall i,j\in A$ $A+i=A+j ,$ and hence ,
$$A+A=\bigcup_{j\in A}(A+j)=A+i .$$
So, $|A+A|=k.$ The claim is proved.

This property is equivalent to the following :
 there exists a point $p\in A$ and a subgroup $H,$ such that $A=p+H.$
So, the points $Vx,V^2x,V^3x,\cdots$ correspond to the sets
$2p+H,4p+H,8p+H,\cdots.$ Actually, such sequences are pre-periodic.

We now  prove the instability of such periodic orbits. Let $l$ be it's
period and $V^l=T$ be the corresponding first return map. Then $T$ has an
instability direction (an eigenvector, whose eigenvalue is greater  than
$1$). It is easy to  check, that  vector $e$ with coordinates 
$$ e_s=\begin{cases}1,\ if\ s\in i+H ;\\
-1\ if \ s\in j+H ;\\
0, ~\mbox{other cases}\\
 \end{cases}
$$
belongs to the plane $\{x:\ \sum x_i=0\}$ (a tangent plane for the
simplex) and realizes this instability direction.

 There is no such directions  for the centre  by the 
 coincidense of each two classes.

So, the basin of attraction for such orbit consists of subvarieties of
strictly positive codimension.We conclude that  almost all orbits tend to the centre of the simplex.

\begin{center}{\bf Remark}\end{center}

We can also consider an infinite-dimensional case. All preceding
arguments remain true. But there is no natural measure on a simplex,and
so we cannot use the term "almost all". 

Another disappointment is  that "most" of the orbits diverge in the 
$l_1$ topology. It follows from the fact, that for "most" orbits the sequence $||V^nx||_{\infty}$ tends to zero. All possible $l_1$ limit points are $l_{\infty}$ limit points. Since the only $l_{\infty}$ limit point is zero, then the subsequence of the   orbit  $l_{1}$ tends to a point, which does not belong to a simplex. This contradiction implies the absense of $l_1$ limit points.

\section*{Acknowledgements}
This research was supported in part by the Uz.R. grant F-2.1.56.


\begin{thebibliography}{99}
\bibitem{1} Bernstein S.N., The solution of a mathematical problem concerning the theory of heredity , Uchenye Zapiski N.-I.Kaf.Ukr.Otd.Mat., {\bf 1}(1924 ) 83-115 (Russian).
\bibitem{2} Lyubich Yu.I., Basic concepts and theorems of the evolution genetics of free populations , Uspekhi Mat.Nauk , {\bf 26}:5 (1971),51-116; English transl. in Russian Math.Surveys, {\bf 26}:5 (1971).
\bibitem{3} Kesten H., Quadratic transformations: a model for population growth.I, II, Adv.Appl.Prob. ,{\bf 2}, 1-82, 179-228 (1970) .
\bibitem{4}Ulam S.M.,{\it A collection of mathematical problems},Interscience Publishers, New York -London (1960) . 
\bibitem{5} Zakharevitch M.I.,On behavior of trajectories and the ergodic hypothesis for quadratic transformations of the simplex, Uspekhi Mat.Nauk ,{\bf 33}(1978 ) 207-208; English transl. in Russian Math.Surveys ,{\bf 33}(1978). 
\bibitem{6} Ganikhodjaev N. N. and Zanin D.V., On ergodicity of quadratic operators on a two dimensional simplex, Uspekhi Mat.Nauk {\bf 59}, no.3(2004)161-162; English transl. in Russian Math.Surveys {\bf 59}(2004) 


\end{thebibliography}
\end{document}